\documentclass[12pt]{article}

\usepackage{amsfonts,amssymb,amsmath,eucal,amsthm}
\usepackage[dvips]{graphics}
\numberwithin{equation}{section}

\newcommand{\Z}{\mathbb{Z}}
\newcommand{\R}{\mathbb{R}}
\newcommand{\C}{\mathbb{C}}
\newcommand{\M}{\mathcal{M}}
\newcommand{\fS}{\mathfrak{S}}

\newcommand{\la}{\lambda}
\newcommand{\de}{\delta}
\newcommand{\sq}{\square}
\newcommand{\sh}{{\textstyle \frac12}}
\newcommand{\da}{\downarrow}

\newcommand{\mc}[2]{\left[#1\to #2\right]}
\newcommand{\F}[4]
  {F\left(\begin{matrix} #1 \,,\, #2 \\ #3 \end{matrix}\, ;
   #4 \right)}
\newcommand{\lan}{\langle}
\newcommand{\ran}{\rangle}
\newcommand{\vac}{\delta_\emptyset}
\newcommand{\LV}{\Lambda^{\frac\infty2}V}
\newcommand{\Y}{\mathbb{Y}}

\newcommand{\slt}{\mathfrak{sl}(2)}
\newcommand{\su}{\mathfrak{su}}

\newtheorem{Theorem}{Theorem}

\theoremstyle{remark}

\begin{document}

\title{$SL(2)$ and $z$-measures}
\author{Andrei Okounkov\thanks{
 Department of Mathematics, University of California at
Berkeley, Evans Hall \#3840, 
Berkeley, CA 94720-3840. E-mail: okounkov@math.berkeley.edu}
}
\date{}

\maketitle

\section{Introduction}

This paper is about the  $z$-measures which are a remarkable
two-parametric family of measures on partitions introduced
in \cite{KOV} in the context of harmonic analysis on the infinite
symmetric group. In a series of papers, A.~Borodin and G.~Olshanski
obtained several fundamental results on these $z$-measures, see their
survey \cite{BO3} which appears in this volume and also \cite{BO4}. The culmination of
this development is an exact determinantal formula for the correlation
functions of the $z$-measures in terms of the hypergeometric
kernel \cite{BO2}. We mention \cite{BOO} as one of the applications
of this formula. The main result of this paper is a representation-theoretic
derivation of the formula of Borodin and Olshanski. 

\subsection{} 

In the early days of $z$-measures, it was already noticed that $z$-measures
have some mysterious connection to representation theory
of $SL(2)$. For example, the $z$-measure is actually positive if
its two parameters $z$ and $z'$ are either complex conjugate $z'=\bar z$ or 
$z,z'\in (n,n+1)$ for some $n\in\Z$. In these cases $z-z'$ is either
imaginary or lies in $(-1,1)$, which was certainly reminiscent of the 
principal and complementary series of representations of $SL(2)$.

Later, S.~Kerov constructed an $SL(2)$-action on partitions for which 
the $z$-measures are certain matrix elements \cite{Ke}.
Finally, Borodin and Olshanski computed the correlation functions
of the $z$-measures is in terms of the Gauss hypergeometric function
which is well known to arise as matrix elements of representations of 
$SL(2)$. The aim of this paper is to put these pieces together.

\subsection{}

I want to thank A.~Borodin, S.~Kerov, G.~Olshanski, and A.~Vershik
for numerous discussions of the $z$-measures. I also want to thank
the organizers of the Random Matrices program at MSRI, especially P.~Bleher,
P.~Deift, and A.~Its. My research was supported by NSF 
under grant DMS-9801466.

The constructions of this paper 
were subsequently  generalized beyond $SL(2)$ and $z$-measures in \cite{O}.

\section{The $z$-measures, Kerov's operators, and correlation functions}\label{s1}

\subsection{Definition of the $z$-measures}

Let $z,z'\in\C$ be two parameters and consider the following measure
on the set of all partitions $\la$ of $n$
\begin{equation}\label{defM}
\M_n(\la)=\frac{n!}{(zz')_n} 
\prod_{\sq\in\la} \frac{(z+c(\sq))(z'+c(\sq))}{h(\sq)^2} \,,
\end{equation}
where 
$$
(x)_n = x(x+1)\dots(x+n-1)\,, 
$$
the product is over all
squares $\sq$ in the diagram of $\la$, $h(\sq)$ is the
length of the corresponding hook, and $c(\sq)$ stands for
the content of the square $\sq$. Recall that, by definition,
the content of $\sq$ is 
$$
c(\sq) = \textup{column}(\sq) -  \textup{row}(\sq)\,,
$$ 
where $\textup{column}(\sq)$ denotes the column number of 
the square $\sq$. The reader is referred to \cite{M}
for general facts about partitions. 

It is not immediately obvious from the definition \eqref{defM} that
\begin{equation}\label{prob}
\sum_{|\la|=n} \M_n(\la)=1 \,.
\end{equation}
One possible proof of \eqref{prob} uses the following 
operators on partitions introduced by S.~Kerov.

\subsection{Kerov's operators}

Consider the vector space with an orthonormal  basis $\{\de_\la\}$
indexed by all partitions of $\la$ of any size. Introduce the
following operators
\begin{alignat*}{2}
U\, \de_\la &= \sum_{\mu=\la+\sq} &(z+c(\sq)) \, &\de_\mu \,\\
L\, \de_\la &= &(zz'+2|\la|)\, &\de_\la\,\\
D\, \de_\la &= \sum_{\mu=\la-\sq} &(z'+c(\sq)) \, &\de_\mu \,,
\end{alignat*}
where $\mu=\la+\sq$ means that $\mu$ is obtained from $\la$ by
adding a square $\sq$ and $c(\sq)$ is the content of this 
square. The letters $U$ and $D$ here stand for ``up'' and ``down''.  

These operators satisfy the commutation relations
\begin{equation}\label{comm}
[D,U]=L\,, \quad [L,U]=2U\,, \quad [L,D]=-2D\,,
\end{equation}
same as for the following basis of $\slt$ 
$$
U=\begin{pmatrix}0 & 1\\ 0 & 0 \end{pmatrix} \,, \quad
L=\begin{pmatrix}1 & 0\\ 0 & -1 \end{pmatrix}\,, \quad
D=\begin{pmatrix}0 & 0\\ -1 & 0 \end{pmatrix}\,.
$$

In particular, it is clear that if $|\la|=n$ then 
$$
(U^n \de_\emptyset, \de_\la) = \dim\la  
\prod_{\sq\in\la} (z+c(\sq))
$$
where
$$
\dim\la = n!  \prod_{\sq\in\la} h(\sq)^{-1} 
$$
is the number of standard tableaux on $\la$. It follows that 
$$
\M_n(\la) = \frac1{n!\, (zz')_n} \,
 (U^n \de_\emptyset,\de_\la) \, (L^n \de_\la, \de_\emptyset)\,. 
$$ 
Using this presentation and the commutation relations  \eqref{comm}
one proves \eqref{prob} by induction on $n$. 

\subsection{The measure $\M$ and its normalization}
In a slightly different language, with induction
on $n$  replaced by the use of generating functions, 
this computation goes
as follows.  

The sequence of the measures $\M_n$ can be conveniently 
assembled as in \cite{BO2} into one measure $\M$ on the set of all partitions
of all numbers as follows
$$
\M = (1-\xi)^{zz'} \sum_{n=0}^\infty \xi^n \, \frac{(zz')_n}{n!} \, \M_n\,,
\quad \xi\in [0,1)\,,
$$
where $\xi$ is a new parameter. In other words, $\M$ is the mixture
of the measures $\M_n$ by means of a negative binomial distribution
on $n$ with parameter $\xi$. 

It is clear that \eqref{prob} 
is now equivalent to $\M$ being a probability measure. It is also clear that
\begin{equation}\label{mM}
\M(\la) =  (1-\xi)^{zz'}  
(e^{\sqrt{\xi}\, U}\, \de_\emptyset,\de_\la) \, (e^{\sqrt{\xi} \, D}\, \de_\la, \de_\emptyset)
\end{equation}
Therefore
\begin{align}\label{Z}
\sum_{\la} \M(\la) = (1-\xi)^{zz'} (e^{\sqrt{\xi} \, D} \, e^{\sqrt{\xi}\, U} \, \de_\emptyset,
\de_\emptyset)
\end{align}

It follows from the definitions that
\begin{equation}\label{Dvac}
D \, \de_\emptyset = 0\,, \quad L \, \de_\emptyset = zz'\, 
\de_\emptyset\,, \quad 
U^* \, \de_\emptyset = 0\,, 
\end{equation}
where $U^*$ is the operator adjoint to $U$. 
Therefore, in order to evaluate \eqref{Z},
it suffices to commute $e^{\sqrt{\xi} \, L}$ through 
$e^{\sqrt{\xi} \, U}$.
 
The following computation in the group $SL(2)$
$$
\begin{pmatrix}
1 & 0 \\
-\beta & 1
\end{pmatrix}
\,
\begin{pmatrix}
1 & \alpha \\
0 & 1
\end{pmatrix}  =
\begin{pmatrix}
1 & \frac{\alpha}{1-\alpha\beta} \\
0 & 1
\end{pmatrix} \,
\begin{pmatrix}
\frac1{1-\alpha\beta} &0 \\
0 & {1-\alpha\beta}
\end{pmatrix} \,
\begin{pmatrix}
1 & 0 \\
-\frac\beta{1-\alpha\beta} & 1
\end{pmatrix}
$$
implies that 
\begin{multline}\label{DU}
\exp\left(\beta\, D\right) \, \exp \left( \alpha\, U \right)=\\ 
\exp\left(\frac\alpha{1-\alpha\beta}\, U \right) \, 
(1-\alpha\beta)^{-L} 
\exp\left(\frac\beta{1-\alpha\beta}\, D\right) \,,
\end{multline}
provided $|\alpha\beta|<1$.  Therefore, 
\begin{align*}
\sum_{\la} \M(\la) &= 
(1-\xi)^{zz'} \left(\exp\left(\frac{\sqrt{\xi}}{1-\xi}\,U\right) \, 
(1-\xi)^{-L} \, \exp\left(\frac{\sqrt{\xi}}{1-\xi}\,D\right)\, \de_\emptyset,
\de_\emptyset\right)\\
&= (1-\xi)^{zz'} \left(
(1-\xi)^{-L} \, \de_\emptyset,
\de_\emptyset\right) = 1 \,,
\end{align*}
as was to be shown. 

\subsection{Correlation functions} 

Introduce the following coordinates
on the set of partitions. To a partition $\la$ we
associate a subset
$$
\fS(\la)=\{\la_i-i+1/2\}\subset \Z+\sh \,.
$$
For example, 
$$
\fS(\emptyset)=\left\{-\frac 12,-\frac 32,-\frac 52,\dots\right\}
$$ 
This set $\fS(\la)$ has the following geometric interpretation. Take the
diagram of $\la$ and rotate it $135^\circ$ as in the following 
picture:
\begin{center}
\scalebox{0.5}{\includegraphics{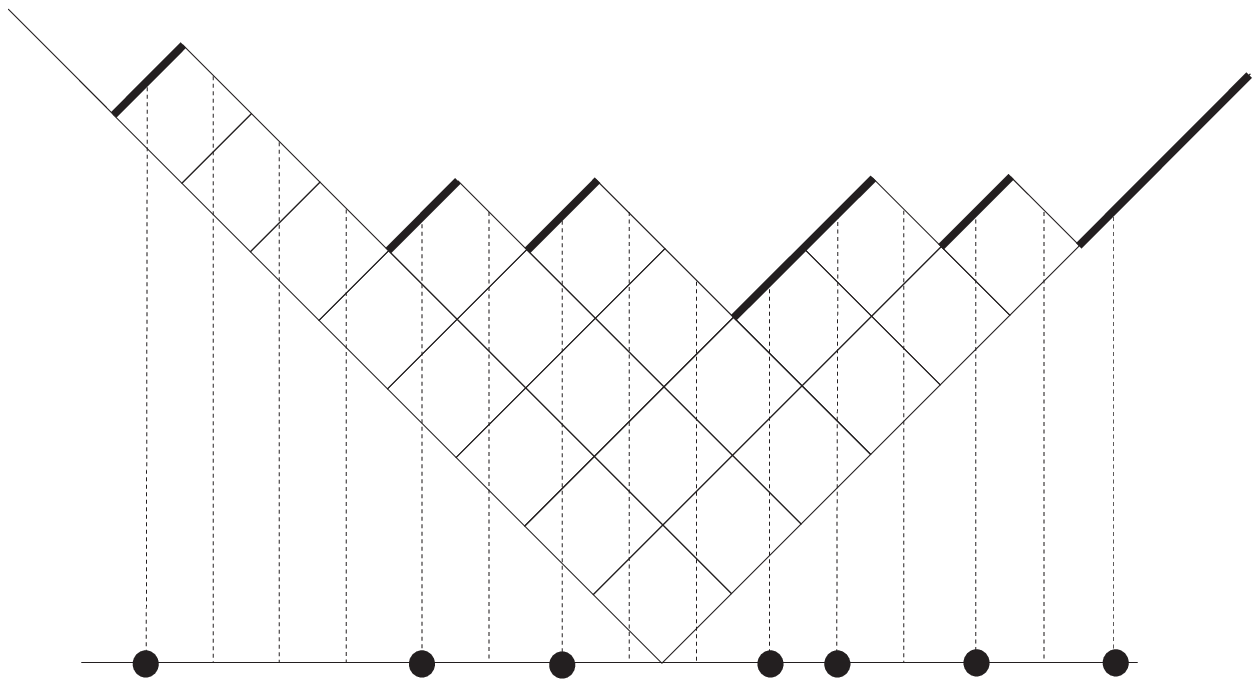}} 
\end{center}
The positive direction of the axis points to the left in the above figure. 
The boundary of $\la$ forms a zigzag path and the elements of $\fS(\la)$,
which  are marked by $\bullet$, correspond to moments when this 
zigzag goes up.

Subsets $S\subset\Z+\frac12$ of the form $S=\fS(\la)$ can be characterized
by
$$
|S_+|=|S_-|<\infty
$$
where
$$
S_+ = S \setminus \left(\Z_{\le 0} - \sh\right) \,, \quad
S_- = \left(\Z_{\le 0} - \sh\right) \setminus S \,. 
$$
The number $|\fS_+(\la)|=|\fS_-(\la)|$ is the number of
squares in the diagonal of the diagram of $\la$ and the
finite set $\fS_+(\la)\cup \fS_-(\la) \subset \Z+\frac 12$
is known as the modified Frobenius coordinates of
$\la$.

Given a finite subset $X\in\Z+\sh$, define the \emph{correlation
function} by 
$$
\rho(X)=\M\big(\{\la, X\subset \fS(\la)\}\big) \,.
$$
In \cite{BO2}, A.~Borodin and G.~Olshanski proved that
$$
\rho(X)=\det\Big[K(x_i,x_j)\Big]_{x_i,x_j\in X} 
$$
where $K$ the \emph{hypergeometric kernel} introduced
in \cite{BO2}. This kernel involves
the Gauss hypergeometric function and the explicit formula for
$K$ will be 
reproduced below.

It is our goal in the present paper
to give a representation-theoretic derivation of 
the formula for correlation functions and, in particular,
show how the kernel $K$ arises from matrix elements of
irreducible $SL(2)$-modules. 

\section{SL(2) and correlation functions} 

\subsection{Matrix elements of $\slt$-modules and Gauss hypergeometric function} 

The fact that the hypergeometric
function arises as matrix coefficients of $SL(2)$ modules is well
known. A standard way to see this is to use a functional
realization of these modules; the computation of 
matrix elements leads then to an integral representation of
the hypergeometric function, see for example how
matrix elements of $SL(2)$-modules are treated in  \cite{Vi}.
An alternative approach is to use explicit formulas
for the action of the Lie algebra $\slt$ and it goes as follows. 

Consider the $\slt$-module $V$ with the basis $v_k$ indexed by
all half-integers $k\in\Z+\sh$ and the following action 
of $\slt$
\begin{alignat}{2}\label{e1d1}
U\, v_k = &&(z+k+\sh)\,&v_{k+1}\,,\\
L\, v_k = &&(2k+z+z')\,&v_{k}\,,\\
D\, v_k = &&(z'+k-\sh)\,&v_{k-1}\,.
\end{alignat}
It is clear that
$$
e^{\alpha\,U}\, v_k = \sum_{s=0}^\infty \frac{\alpha^s}{s!}\,
(z+k+\sh)_s \, v_{k+s} \,.
$$
Introduce the following notation
$$
(a)_{\da s} = a(a-1)(a-2)\cdots (a-s+1) \,.
$$
With this notation we have
$$
e^{\beta\, D} \, v_k = \sum_{s=0}^\infty \frac{\beta^s}{s!}\,
(z'+k-\sh)_{\da s} \, v_{k-s} \,.
$$
Denote by $\mc{i}{j}_{\alpha,\beta,z,z'}$ the coefficient
of $v_j$ in the expansion of $e^{\alpha \,U} \, e^{\beta\, D}\, v_i$
$$
e^{\alpha \,U} \, e^{\beta\, D}\, v_i = \sum_j \,
\mc{i}{j}_{\alpha,\beta,z,z'}
\, v_j \,.
$$
A direct computation yields  
\begin{multline} \label{e12} 
\mc{i}{j}_{\alpha,\beta,z,z'}= \\ 
\begin{cases}
{\displaystyle \frac{\alpha^{j-i}}{(j-i)!} \, (z+i+\sh)_{j-i}} \, 
\F{-z-i+\sh}{-z'-i+\sh}{j-i+1}{\alpha\beta} \,, & i \le j  \,, \\
{\displaystyle \frac{\beta^{i-j}}{(i-j)!} \, (z'+j+\sh)_{i-j}} \, 
\F{-z-j+\sh}{-z'-j+\sh}{i-j+1}{\alpha\beta} \,, 
& i \ge j  \,,
\end{cases}
\end{multline}
where 
$$
\F{a}{b}{c}{z}=\sum_{k=0}^\infty \frac{(a)_k\, (b)_k}{(c)_k \, k!} \,
z^k \,,
$$
is the Gauss hypergeometric function.

Consider now the dual module $V^*$ spanned by functionals $v^*_j$ such
that
$$
\lan v^*_i, v_j \ran=\de_{ij}
$$ 
and equipped with the dual action of $\slt$
\begin{alignat*}{2}
U\, v^*_k &= - (z+k-\sh)&&\,v^*_{k-1} \,, \\
D\, v^*_k &= - (z'+k+\sh)&&\,v^*_{k+1} \,.
\end{alignat*}
Denote by $\mc{i}{j}^*_{\alpha,\beta,z,z'}$ the coefficient
of $v^*_j$ in the expansion of $e^{\alpha \,U} \, e^{\beta\, D}\, v^*_i$
$$
e^{\alpha \,U} \, e^{\beta\, D}\, v^*_i = \sum_j \,
\mc{i}{j}^*_{\alpha,\beta,z,z'}
\, v^*_j \,.
$$
We have
\begin{multline} \label{e13} 
\mc{i}{j}^*_{\alpha,\beta,z,z'}= \\ 
\begin{cases}
{\displaystyle \frac{(-\beta)^{j-i}}{(j-i)!} \, (z'+i+\sh)_{j-i}} \, 
\F{z+j+\sh}{z'+j+\sh}{j-i+1}{\alpha\beta} \,, & i \le j \,, \\
{\displaystyle \frac{(-\alpha)^{i-j}}{(i-j)!} \, (z+j+\sh)_{i-j}} \, 
\F{z+i+\sh}{z'+i+\sh}{i-j+1}{\alpha\beta} \,, 
& i \ge j \,,
\end{cases}
\end{multline}

\subsection{Remarks}
\subsubsection{Periodicity} 
Observe that representations whose parameters $z$ and $z'$
are related by the tranformation
$$
(z,z')\mapsto (z+m,z'+m)\,, \quad m\in\Z\,, 
$$
are equivalent. The above transformation amounts to just
a renumeration of the vectors $v_k$. G.~Olshanski  pointed
out that this periodicity in $(z,z')$ is reflected in a
similar periodicity of various asymptotic properties of
$z$-measures, see Sections 10 and 11 of \cite{BO4}. 

\subsubsection{Unitarity}

Recall that the $z$-measures are  positive if either
$z'=\bar z$ or $z,z'\in(n,n+1)$ for some $n$. By analogy
with representation theory of $SL(2)$, these cases were
called the principal and the complementary series. 

Observe that in these case the above representations 
have a positive defined Hermitian form $Q$ which is
invariant in the following sense
$$
Q(Lu,v)=Q(u,Lv)\,, \quad Q(Uu,v)=Q(u,Dv)\,.
$$
The form $Q$ is given by
$$
Q(v_k,v_k)=
\begin{cases}1  & z'=\bar z\,, \\
\dfrac{\Gamma(z'+k+\frac12)}{\Gamma(z+k+\frac12)}\, 
 &z,z'\in(n,n+1)\,,
     \end{cases}
$$
and $Q(v_k,v_l)=0$ if $k\ne l$. It follows that the operators
$$
\tfrac{i}2\,L,\tfrac12\,(U-D), \tfrac{i}2\,(U+D) \in \slt\,,
$$
which form a standard basis of $\su(1,1)$, are skew-Hermitian
and hence this representation of $\su(1,1)$ can be integrated
to a unitary representation of the universal covering
group of $SU(1,1)$. This group $SU(1,1)$ is isomorphic to
$SL(2,\R)$ and the above representations correspond to
the principal and complementary series of unitary
representations of the universal covering of
$SL(2,\R)$, see \cite{Pu}.

\subsection{The infinite wedge module}

Consider the module $\Lambda^{\frac{\infty}2}\, V$ which is, by
definition, spanned by vectors
$$
\delta_S=v_{s_1} \wedge v_{s_2} \wedge  v_{s_3} \wedge  \dots\,,
$$
where $S=\{s_1>s_2>\dots\}\subset \Z+\sh$ is a such subset that
both sets
$$
S_+ = S \setminus \left(\Z_{\le 0} - \sh\right) \,, \quad
S_- = \left(\Z_{\le 0} - \sh\right) \setminus S 
$$
are finite. We equip this module with the inner product 
in which the basis $\{\delta_S\}$ is orthonormal. Introduce the
following operators 
$$
\psi_k, \psi^*_k : \Lambda^{\frac{\infty}2}\, V \to
\Lambda^{\frac{\infty}2}\, V \,.
$$
The operator $\psi_k$ is the exterior multiplication by $v_k$
$$
\psi_k \left(f\right) = v_k \wedge f  \,. 
$$
The operator $\psi^*_k$ is the adjoint operator; it can be
also given by the formula
$$
\psi^*_k \left(v_{s_1} \wedge v_{s_2} \wedge  v_{s_3} \right)=
\sum_i (-1)^{i+1} \lan v^*_k, v_{s_i} \ran \,
v_{s_1} \wedge v_{s_2} \wedge \dots \wedge \widehat{v_{s_i}} \wedge
\dots \,.
$$
These operators satisfy the canonical anticommutation relations
$$
\psi_k \psi^*_k + \psi^*_k \psi_k = 1\,,
$$
all other anticommutators being equal to $0$. It is clear that
\begin{equation}\label{e14}
\psi_k \psi^*_k \, \delta_S = 
\begin{cases}
\delta_S\,, & k \in S \,, \\
0 \,, & k \notin S \,.
\end{cases}
\end{equation}
A general reference on the infinite wedge space 
is Chapter 14 of the book \cite{K}. 

The Lie algebra $\slt$ acts on $\Lambda^{\frac{\infty}2}\, V$.
The action of $U$ and $D$ are the obvious extensions of the action
on $V$. In terms of the fermionic operators $\psi_k$ and $\psi^*_k$
they can be written as follows
\begin{align*}
U&=\sum_{k\in\Z+\frac12} (z+k+\sh)\,\, \psi_{k+1} \psi^*_{k} \,,\\
D&=\sum_{k\in\Z+\frac12} (z'+k+\sh)\,\, \psi_{k} \psi^*_{k+1} \,.
\end{align*}
The easiest way to define the action of $L$ is to set it equal to
$[D,U]$ by definition. We obtain
$$
L=2H + (z+z')\, C + zz' \,,
$$
where $H$ is the energy operator
$$
H=\sum_{k>0} k \, \psi_k \psi^*_k - \sum_{k<0} k \, \psi^*_k \psi_k
\,,
$$
and $C$ is the charge
$$
C=\sum_{k>0} \psi_k \psi^*_k - \sum_{k<0}  \psi^*_k \psi_k \,.
$$
It is clear that
$$
C\, \delta_S = \left(|S_+|-|S_-|\right) \, \delta_S
$$
and, similarly,
$$
H \, \delta_S =  \left(\sum_{k\in S_+} k - \sum_{k \in S_-} k\right)\,
\delta_S\,.
$$
The charge is preserved by the $\slt$ action. 

Consider the zero charge subspace, that is, the kernel of $C$
$$
\Lambda_0 \subset \LV\,.
$$
It is spanned by vectors which, abusing notation, we shall
denote by
\begin{equation}\label{e15}
\delta_\la=\delta_{S(\la)} \,, \quad S(\la)=
\left\{\la_1-\tfrac12,\la_2-\tfrac32,\la_3-\tfrac52,\dots \right\}\,,
\end{equation}
where $\la$ is a partition. One immediately sees that
the action of $\slt$ on $\{\delta_\la\}$ is identical with 
Kerov's operators. 

\subsection{Correlation functions}

Recall that the correlation functions were defined by
$$
\rho(X)=\M\big(\{\la, X\subset \fS(\la)\}\big) \,,
$$
where the finite set
$$
X=\{x_1,\dots,x_s\} \subset \Z+\tfrac12 
$$
is arbitrary. 

The important observation is that
\eqref{mM} and \eqref{e14}  imply the following
expression for the correlation functions 
\begin{equation}\label{e17}
\rho(X) =  (1-\xi)^{zz'} \, 
\left( e^{\sqrt\xi\, D} \, \prod_{x\in X} \psi_{x} \psi^*_{x} \,
 e^{\sqrt\xi\, U}\, \vac,\vac\right) \,.
\end{equation}
We apply to \eqref{e17} the same strategy we applied to \eqref{Z} 
which is to commute the operators $e^{\sqrt\xi\, D}$ and $e^{\sqrt\xi\, U}$ 
all the way to the 
right and left, respectively, and then use \eqref{Dvac}. 
{F}rom \eqref{DU}, we have for any operator $A$ the following identity
\begin{multline}\label{e18}
e^{\beta\, D}\, A \,e^{\alpha\, U}= \\
e^{\frac\alpha{1-\alpha\beta}\, U} \,
\left[
e^{-\frac\alpha{1-\alpha\beta}\, U}\,
e^{\beta\, D} \, A \,
e^{-\beta\, D}
e^{\frac\alpha{1-\alpha\beta}\, U}
\right]
(1-\alpha\beta)^{-L} e^{\frac\beta{1-\alpha\beta}\, D} \,.
\end{multline}
We now apply this identity with $\alpha=\beta=\sqrt\xi$
and $A=\prod \psi_x \psi^*_x$ to obtain 
\begin{equation}\label{e19}
\rho(X) =  
\left(G\, \prod_{x\in X} \psi_{x} \psi^*_{x}  \, G^{-1}\, \vac,\vac\right) \,,
\end{equation}
where
$$
G=\exp\left(\frac{\sqrt{\xi}}{\xi-1}\, U\right)\,
\exp\left(\sqrt{\xi}\, D\right)\,.
$$

Consider the following operators
\begin{alignat}2
\Psi_k &= G\, \psi_k \, G^{-1} &&= \sum_i
\mc{k}{i}\, \psi_i \,, \label{e110}\\
\Psi^*_k &= G\, \psi^*_k \, G^{-1} &&= \sum_i
\mc{k}{i}^*\, \psi^*_i \,, \label{e111} 
\end{alignat}
with the understanding that matrix elements without 
parameters stand for the following choice of parameters
\begin{equation}\label{e112}
\mc{k}{i}=\mc{k}{i}_{\xi^{1/2}(\xi-1)^{-1},\, \xi^{1/2},z,z'}\,,
\end{equation}
and with same choice of parameters for $\mc{k}{i}^*$.
The first equality in both \eqref{e110} and \eqref{e111} is a definition
and the second equality follows from the definition of the operators
$\psi_i$ and the definition of the 
matrix coefficients $\mc{i}{j}_{\alpha,\beta,z,z'}$.

From \eqref{e19} we obtain
\begin{equation}\label{e114}
\rho(X)=\left(\prod_{x\in X} \Psi_{x} \Psi^*_{x}  \, \vac,\vac\right) \,.
\end{equation}
Applying Wick's theorem to \eqref{e114}, or 
simply unraveling the definitions in the right-hand side
of \eqref{e114}, we obtain the following

\begin{Theorem} We have 
\begin{equation}
\rho(X)=\det\big[K(x_i,x_j)\big]_{1\le i,j \le s} \,,
\end{equation}
where the kernel $K$ is defined by 
$$
K(i,j)=\left(\Psi_{i} \Psi^*_{j} \, \vac,\vac\right) \,.
$$
\end{Theorem}

Observe that
$$
\left(\psi_l \, \psi^*_m \, \vac,\vac\right)=
\begin{cases}
1\,, & l=m < 0 \,, \\
0\,, & \text{otherwise} \,.
\end{cases}
$$
Therefore, applying the formulas \eqref{e110} and \eqref{e111} we obtain

\begin{Theorem}  We have 
\begin{equation}\label{e115}
K(i,j)=\sum_{m=-1/2,-3/2,\dots}
\mc{i}{m} \, \mc{j}{m}^* \,,
\end{equation}
with the agreement \eqref{e112} about  matrix elements without 
parameters.
\end{Theorem}

 The formula \eqref{e115} is the analog of the 
Proposition 2.9  in \cite{BOO} for the discrete Bessel kernel. 

We conclude this section with the following formula which,
after substituting the formulas \eqref{e12} and \eqref{e13} for
matrix elements, becomes the formula of Borodin and Olshanski \cite{BO2}. 

\begin{Theorem} We have
\begin{equation}\label{e116}
K(i,j)=\frac{
z' \sqrt{\xi} \mc{i}{\sh}\, \mc{j}{-\sh}^*  - 
z \frac{\sqrt{\xi}}{(\xi-1)^2} \mc{i}{-\sh}\, \mc{j}{\sh}^* }{i-j} \,,
\end{equation}
where for $i=j$ the right-hand side is defined 
by continuity. 
\end{Theorem}

More generally, set
$$
K(i,j)_{\alpha,\beta} = \left(\Psi_i(\alpha,\beta)\,
\Psi^*_j(\alpha,\beta) \, \vac,\vac \right)\,.
$$
where
\begin{alignat*}2
\Psi_k &= e^{\alpha\, U} \, e^{\beta\, D}\, \psi_k \, 
e^{-\beta\, D} \, e^{-\alpha\, U} &&= \sum_i
\mc{k}{i}_{\alpha,\beta,z,z'}\, \psi_i \\
\Psi^*_k &= e^{\alpha\, U} \, e^{\beta\, D}\, \psi^*_k \, 
e^{-\beta\, D} \, e^{-\alpha\, U} &&= \sum_i
\mc{k}{i}^*_{\alpha,\beta,z,z'}\, \psi^*_i \,. 
\end{alignat*}

We will prove that 
\begin{multline}\label{e117} 
K(i,j)_{\alpha,\beta}=
\left(
\beta z' \mc{i}{\sh}_{\alpha,\beta,z,z'}\, \mc{j}{-\sh}^*_{\alpha,\beta,z,z'}  -
\right.\\
\left.
\alpha (\alpha\beta-1) z 
 \mc{i}{-\sh}_{\alpha,\beta,z,z'}\, \mc{j}{\sh}^*_{\alpha,\beta,z,z'} \right)
/(i-j) \,.
\end{multline}

First, we treat the case $i\ne j$ in which case we can
clear the denominators in \eqref{e117}. From the following
computation with $2\times 2$ matrices
\begin{multline*}
\begin{pmatrix}
1 & \alpha \\
0 & 1
\end{pmatrix} 
\begin{pmatrix}
1 & 0 \\
-\beta & 1
\end{pmatrix} 
\begin{pmatrix}
1 & 0 \\
0 & -1
\end{pmatrix} 
\begin{pmatrix}
1 & 0 \\
\beta & 1
\end{pmatrix} 
\begin{pmatrix}
1 & -\alpha \\
0 & 1
\end{pmatrix} = \\
\begin{pmatrix}
1-2\alpha\beta & 2\alpha(\alpha\beta-1) \\
-2\beta & 2\alpha\beta-1
\end{pmatrix}
\end{multline*}
we conclude that
$$
e^{\alpha\, U} \, e^{\beta\, D} \, L 
e^{-\beta\, D} \, e^{-\alpha\, U} 
=L+T\,,
$$
where
$$
T=-2\alpha\beta\, L + 2\beta\, D +2\alpha(\alpha\beta-1) \, U \,.
$$
This can be rewritten as follows
\begin{align}
[L,e^{\alpha\, U} \, e^{\beta\, D}]&= -T\,
e^{\alpha\, U} \, e^{\beta\, D}\,, \label{e118}\\
[L,e^{-\alpha\, U} \, e^{-\beta\, D}]&=
e^{-\alpha\, U} \, e^{-\beta\, D} \, T\,. \label{e119}
\end{align}
Note that
\begin{equation}\label{e120}
[L,\psi_i \, \psi^*_j] = 2(i-j)\, \psi_i \, \psi^*_j \,.
\end{equation}
From \eqref{e118}, \eqref{e119}, and \eqref{e120} we have 
\begin{multline}\label{e121}
\left[L,\Psi_i(\alpha,\beta)\,
\Psi^*_j(\alpha,\beta) \right]= \\
-\left[T,\Psi_i(\alpha,\beta)\,
\Psi^*_j(\alpha,\beta) \right]+2(i-j) \, \Psi_i(\alpha,\beta)\,
\Psi^*_j(\alpha,\beta) \,.
\end{multline}

Since $\vac$ is an eigenvector of $L$ we have
$$
\left(\left[L,\Psi_i(\alpha,\beta)\,
\Psi^*_j(\alpha,\beta) \right]\, \vac,\vac \right) = 0
$$
Expand this equality using \eqref{e121} and the relations
\begin{align*}
T \, 
\vac &= - 2\alpha\beta z z' \vac + 2\alpha(\alpha\beta-1) z \delta_\sq\,,\\
T^*
\vac &=  - 2\alpha\beta z z' \vac+ 2 \beta z' \delta_\sq \,,
\end{align*}
where $T^*$ is the operator adjoint to $T$ and $\delta_\sq$ is the 
vector corresponding to the partition $(1,0,0,\dots)$.
We obtain
\begin{multline*}
(i-j) K(i,j)_{\alpha,\beta} =
 \beta z'  \, \left(\Psi_i(\alpha,\beta)\,
\Psi^*_j(\alpha,\beta)\, \vac,\delta_\sq \right) - \\
\alpha(\alpha\beta-1) z \,
\left(\Psi_i(\alpha,\beta)\,
\Psi^*_j(\alpha,\beta)\, \delta_\sq,\vac \right)
\end{multline*}
In order to obtain \eqref{e117} for $i\ne j$, it now
remains to observe that
\begin{align*}
(\psi_l\, \psi^*_m \, \vac,\delta_\sq)&=
\begin{cases}
1\,, & l=\frac12\,, \,\, m =-\frac12 \,,\\
0 \,, & \text{otherwise}\,,
\end{cases}
\\ 
(\psi_l\, \psi^*_m \, \delta_\sq,\vac)&=
\begin{cases}
1\,, & l=-\frac12\,, \,\,  m =\frac12 \,,\\
0 \,, & \text{otherwise}\,.
\end{cases}
\end{align*}

In the case $i=j$ we argue by continuity. It is clear
from \eqref{e115} that $K(i,j)$ is an analytic function
of $i$ and $j$ and so is the right-hand side of \eqref{e116}.
The passage from \eqref{e115} to \eqref{e116} is based on the 
fact that the product $i$ times $\mc{i}{m}_{\alpha,\beta,z,z'}$ is 
a linear combination of $\mc{i}{m}_{\alpha,\beta,z,z'}$ and
$\mc{i}{m\pm 1}_{\alpha,\beta,z,z'}$ with coefficients
which are linear functions of $m$. Since the matrix
coefficients are, essentially, the hypergeometric
function, such a relation must hold for any $i$, not just
half-integers. Hence, \eqref{e115} and \eqref{e116} are equal
for any $i\ne j$, not necessarily half-integers. 
Therefore, they are equal for $i=j$.

\subsection{Rim-hook analogs}

The same principles
apply to rim-hook analogs of the $z$-measures
which were also considered by S.~Kerov \cite{Ke}.

Recall that a rim hook of a diagram $\la$ is, by
definition, a skew diagram $\la/\mu$ which is connected 
and lies on the rim of $\la$.  Here connected
means that the squares have to be connected by
common edges, not just common vertices. Rim hooks
of a diagram $\la$ are in the following 1-1 correspondence with 
the squares of $\la$: given
a square $\sq\in\la$, the corresponding rim hook
consists of all squares on the rim of $\la$ which
are (weakly) to the right of and below $\sq$.
The length of this rim hook is equal to the
hook-length of $\sq$.  

The entire discussion of the previous section applies to the 
more general operators
\begin{alignat*}{2}
U_r\, v_k = &&\left(z+\tfrac kr +\sh\right)\,&v_{k+r}\,,\\
L_r\, v_k = &&\left(\tfrac{2k}r+z+z'\right)\,&v_{k}\,,\\
D_r\, v_k = &&\left(z'+\tfrac kr-\sh\right)\,&v_{k-r}\,,
\end{alignat*}
which satisfy the same $\slt$ commutation relations. The
easiest way to check the commutation relations
 is to consider  $\frac kr$ rather
than $k$ as the index of $v_k$; the above formulas then
become precisely the formulas \eqref{e1d1}. The 
operator $U_r$ acts on the basis $\{\delta_\la\}$ as follows 
\begin{equation}\label{e401}
U_r \, \delta_\la=\sum_{\mu=\la+\textup{ rim hook}} (-1)^{\textup{height}+1} 
\left(z+\frac1{r^2}\sum_{\sq\in\textup{rim hook}}
c(\sq) \right) \delta_\mu \,,
\end{equation}
where the summation is over all partitions $\mu$ which
can be obtained from $\la$ by adding a rim hook of length
$r$, height is the number of horizontal rows occupied
by this rim hook and $c(\sq)$ stands, as usual for the
content of the square $\sq$. Similarly, the operator
$D_r$ removes rim hooks of length $r$. These operators
were considered by Kerov \cite{Ke}. 

It is clear that the action of the operators 
$e^{\alpha U_r}$ and $e^{\beta D_r}$
on a half-infinite wedge product like
$$
v_{s_1} \wedge v_{s_2} \wedge  v_{s_3} \wedge  \dots\,,
$$
essentially (up to a sign which disappears in formulas
like \eqref{e17}) factors into the tensor
product of $r$ separate actions on
$$
\bigwedge_{s_i\equiv k+\frac12\!\! \mod r} v_{s_i} \,, \quad k=0,\dots,r-1 \,.
$$
Consequently, the analogs of the correlation functions
\eqref{e17} have again a determinantal form with a
certain kernel $K_r(i,j)$ which has the following structure.
If $i\equiv j \mod r$ then $K_r(i,j)$ is essentially the
 kernel $K(i,j)$ with rescaled arguments.
Otherwise, $K_r(i,j)=0$. 
 
This factorization of the action on $\LV$ is just one
more way to understand the following well-known 
phenomenon. Let $\Y_r$ be the partial ordered set formed
by partitions with respect to the following ordering:
$\mu\le_r\la$ if $\mu$ can be obtained from $\la$ by
removing a number of rim hooks with $r$ squares. The
minimal elements of $\Y_r$ are called the $r$-cores.
The $r$-cores are precisely those partitions which
do not have any hooks of length $r$. We have 
\begin{equation}\label{e41}
\Y_r \cong \bigsqcup_{\textup{$r$-cores}} (\Y_1)^r  
\end{equation}
as partially ordered sets. Here the
Cartesian product $(\Y_1)^r$ is ordered as follows:
$$
(\mu_1,\dots,\mu_r) \le (\la_1,\dots,\la_r)
\quad  \Leftrightarrow \quad \mu_i \le_1 \la_i, \quad 
i=1,\dots,r \,,
$$
and the partitions corresponding to different
$r$-cores are incomparable in the $\le_r$-order. 
Combinatorial algorithms
which materialize the isomorphism \eqref{e41}
are discussed in Section 2.7 of the book \cite{JK}.
The $r$-core and the $r$-tuple of partitions which
the isomorphism \eqref{e41} associates to a partition
$\la$ are called the $r$-core of $\la$ and the 
$r$-quotient of $\la$. Among more recent papers dealing with $r$-quotients
let us mention \cite{FS}
where an approach similar to the use of $\LV$
is employed, an analog of the Robinson-Schensted
algorithm for $\Y_r$ is discussed, and further
references are given. 

The factorization \eqref{e41} and the corresponding
analog of the Robinson-Schensted
algorithm play the central role in the recent paper
\cite{B}, see also \cite{R}.

\end{document}